\documentclass{elsart1}

 \usepackage{graphicx}

\usepackage{amssymb}

\newcommand{\R}{{\Bbb R}}

\begin{document}
\begin{frontmatter}
\title{Monotone traveling wavefronts of the KPP-Fisher delayed equation}
\author[c]{Adrian Gomez}
\author[c]{and Sergei Trofimchuk}
\address[c]{Instituto de Matem\'atica y Fisica, Universidad de Talca, Casilla 747,
Talca, Chile \\ {\rm adriangomez79@hotmail.com and
trofimch@imath.kiev.ua}}

\bigskip

\begin{abstract}
\noindent In the early 2000's, Gourley (2000),  Wu {\it et al.}
(2001), Ashwin {\it et al.} (2002) initiated the study of the
positive wavefronts in the delayed
Kolmogorov-Petrovskii-Piskunov-Fisher equation
$$u_t(t,x) = \Delta u(t,x)  + u(t,x)(1-u(t-h,x)), \ u \geq 0,\ x
\in \R^m. \eqno(*)$$ Since then, this model has become one of the
most popular objects in the studies of  traveling waves for the
monostable delayed reaction-diffusion equations. In this paper, we
give a complete solution to the problem of existence and uniqueness
of monotone waves in equation $(*)$. We show that each monotone
traveling wave can be found via an iteration procedure. The proposed
approach is based on the use of special monotone integral operators
(which are different from the usual Wu-Zou operator) and appropriate
upper and lower solutions associated to them. The analysis of the
asymptotic expansions of the eventual traveling fronts at infinity
is another key ingredient of our approach.
\end{abstract}
\end{frontmatter}
\newpage

\section{Introduction and main results}\vspace{-3mm}
\label{intro} It is well known that the traveling waves theory was
initiated in 1937  by Kolmogorov, Petrovskii, Piskunov \cite{KPP}
and Fisher \cite{Fish} who studied the wavefront solutions of the
diffusive logistic equation
\begin{equation}\label{kppf} \hspace{5mm}
u_t(t,x) = \Delta u(t,x)  + u(t,x)(1-u(t,x)), \ u \geq 0,\ x \in
\R^m. \end{equation} \vspace{-5mm} We recall that the classical
solution $u(x,t) = \phi(\nu \cdot  x +ct),\ \|\nu\| =1,$ is a
wavefront (or a traveling front) for (\ref{kppf}), if the profile
function $\phi$ is positive and satisfies $\phi(-\infty) = 0, \
\phi(+\infty) = 1$.

The existence of the wavefronts in (\ref{kppf}) is equivalent to the
presence of  positive heteroclinic connections in an associated
second order non-linear differential equation. The phase plane
analysis is the natural geometric way to study these heteroclinics.
The method is conclusive enough to demonstrate that (a) for every $c
\geq 2$, the KPP-Fisher equation has exactly one traveling front
$u(x,t) = \phi(\nu \cdot x +ct)$; (b) Eq. (\ref{kppf}) does not have
any traveling front propagating at the velocity $c < 2$; (c) the
profile $\phi$ is necessarily strictly increasing function.

The stability of traveling fronts in (\ref{kppf}) represents another
important aspect of the topic: however, we do not discuss it here.
Further reading and relevant information can be found in
\cite{brams,Lau,RK,yana}.

Eq. (\ref{kppf}) can be viewed as a natural extension of the
ordinary logistic equation $u'(t) = u(t)(1-u(t))$. An important
improvement of this growth model was proposed by Hutchinson
\cite{hut} in 1948 who incorporated the maturation delay $h > 0$ in
the following way:
\begin{equation}\label{logh} \hspace{5mm}
u'(t) = u(t)(1-u(t-h)), \ u \geq 0.
\end{equation}
\vspace{-3mm} This  model is now commonly known as the Hutchinson' s
equation. Since then, the delayed KPP-Fisher equation or the
diffusive Hutchinson's equation
\begin{equation}\label{17} \hspace{5mm}
u_t(t,x) = \Delta u(t,x)  + u(t,x)(1-u(t-h,x)), \ u \geq 0,\ x \in
\R^m, \end{equation} \vspace{-5mm} is considered as a natural
prototype of delayed reaction-diffusion equations. It has attracted
the attention of many authors, see
\cite{Ai,ZAMP,fhw,Fri,gop,hadeler,lu,wlr,wz,y}. In particular, the
existence of  traveling fronts connecting the trivial and positive
steady states in (\ref{17}) (and its non-local generalizations) was
studied in \cite{Ai,ZAMP,cdm,fhw,gou2000,pan,wlr,wz}. Observe that
the biological meaning of $u$ is the size of an adult population,
therefore {\it only} non-negative solutions of (\ref{17}) are of
interest. It is worth to mention that there is another delayed
version of Eq. (\ref{kppf}) derived by  Kobayashi \cite{Koba} from a
branching process:
$$
u_t(t,x) = \Delta u(t,x)  + u(t-h,x)(1-u(t,x)), \ u \geq 0,\ x \in
\R^m.
$$
However, since the right-hand side of this equation is monotone
increasing with respect to the delayed term, the theory of this
equation is fairly different (and seems to be simpler) from the
theory of (\ref{17}), see \cite{sch,wz,z}.

This paper deals with the problem of existence and uniqueness of
monotone  wavefronts for Eq. (\ref{17}). The phase plane analysis
does not work now because of the infinite dimension of phase spaces
associated  to delay equations. Recently, the existence problem was
considered by using two different approaches. The first method,
which was proposed in \cite{wz}, uses the positivity and
monotonicity properties of the integral operator
\begin{equation}\label{psea} \hspace{-3mm}
(A\phi)(t) = \frac{1}{\epsilon'}\left\{\int_{-\infty}^te^{r_1
(t-s)}(\mathcal{H}\phi)(s)ds + \int_t^{+\infty}e^{r_2
(t-s)}(\mathcal{H}\phi)(s)ds \right\},
\end{equation}
\vspace{-5mm} where $(\mathcal{H}\phi)(s)= \phi(s)(\beta+1
-\phi(s-h))$ for some appropriate $\beta >1$, and $\epsilon' =
\epsilon(r_2-r_1)$ with $r_1<0<r_2$ satisfying $\epsilon z^2 -z
-\beta =0$, and $\epsilon^{-1/2} = c
> 0$ is the front velocity. A direct verification shows that the
profiles $\phi \in C(\R,\R_+)$ of traveling waves are completely
determined by the integral equation $A\phi=\phi$. Wu and Zou have
found a subtle combination of the usual and the Smith and Thieme
nonstandard orderings on an appropriate profile set $\Gamma^*
\subset C(\R,(0,1))$ which allowed them (under specific
quasimonotonicity conditions) to indicate a pair of upper and lower
solutions $\phi^\pm$ such that $\phi^-\leq A^{j+1}\phi^+ \leq
A^{j}\phi^+, \ j=0,1, \dots$ Then the required traveling front
profile is given by $\phi = \lim A^{j}\phi^+$. More precisely, in
\cite[Theorem 5.1.5]{wz}, Wu and Zou established the following
\begin{prop} \label{1wz}
For any $c>2$, there exists $h^*(c)>0$  such that if $h \leq
h^*(c)$, then Eq. (\ref{17}) has a monotone traveling front with
wave speed $c$.
\end{prop}
The above result was  complemented in \cite[Remark 5.15]{wlr} and
\cite{pan}, where it was shown that Proposition \ref{1wz} remains
valid if $c=2$. It should be observed that Wang {\it et al.}
\cite{wlr} have also used the method of upper and lower solutions,
however their lower solution is different from that in \cite{wz}.
Recently, Ou and Wu \cite{OW} showed that Proposition \ref{1wz} can
be proved by means of a perturbation argument (considering $h>0$ as
a small parameter). 

The second method was proposed in \cite{fhw}. It essentially relies
on the fact that, in a 'good' Banach space, the Frechet derivative
of $\lim_{\epsilon \to 0} A$ along a heteroclinic solution $\psi$ of
the limit delay differential equation (\ref{logh}) is a surjective
Fredholm operator. In consequence, the Lyapunov-Schmidt reduction
was used to prove the existence of a smooth family of wave solutions
in some neighborhood of $\psi$. The following result was proved in
\cite[Corollary 6.6.]{fhw}:
\begin{prop} \label{2wz}
There exists $c^* >0$ such that if $0<h<1/e$ then for any $c>c^*$,
Eq. (\ref{17}) has a wave solution $u(x,t) = \phi(\nu \cdot  x
+ct),\ |\nu| =1,$ satisfying $\phi(-\infty) = 0, \ \phi(+\infty) =
1$.
\end{prop}

We remark that the positivity of this wave was not proved in
\cite{fhw} and the value of $c^*>0$ was not given explicitly.
Nevertheless, as it was shown in \cite{FT} for the case of the
Mackey-Glass type equations, the method of \cite{fhw} may be refined
to establish the existence of positive wavefronts as well. Moreover,
it follows from \cite{FT} that Proposition \ref{2wz} is still valid
for $h \in (0,3/2)$. The recent work \cite{AVT} suggests that the
approach of \cite{fhw} can be also used to prove the uniqueness (up
to shifts) of the positive traveling solution of (\ref{17}) for
sufficiently fast speeds.

In this paper, motivated by ideas in \cite{DK,wz}, we give a
criterion for the existence of positive monotone wavefronts in
(\ref{17}) and prove their uniqueness (modulo translation). In order
to do this, instead of using operator (\ref{psea}) as it was done in
all previous works, we work with different integral operators,
namely:
\begin{equation}\label{iie} \hspace{5mm}
(\mathcal{A}\varphi)(t) = \frac{1}{\epsilon(\mu-\lambda)}
\int_t^{+\infty}(e^{\lambda (t-s)}- e^{\mu
(t-s)})\varphi(s)\varphi(s-h)ds,
\end{equation}
\vspace{-5mm} where $\epsilon \in (0,0.25)$ and $0 < \lambda < \mu$
are the roots of $\epsilon z^2 -z + 1 =0$, and with
\begin{equation}\label{iie2} \hspace{5mm}
(\mathcal{B}\varphi)(t) = 4\int_t^{+\infty}(s-t)e^{2
(t-s)}\varphi(s)\varphi(s-h)ds
\end{equation}
\vspace{-2mm} which can be considered as the limit of $\mathcal{A}$
when $\epsilon \to 0.25$.
 Remarkably,  all monotone
wavefronts (in particular, the wavefronts propagating with the
minimal speed $c=2$) can be found via a monotone iterative algorithm
which uses $\mathcal{A}, \mathcal{B}$ and converges uniformly on
$\R$.

Before stating our main results, let us introduce the critical delay
$h_1=  0.560771160\dots$ This value coincides with the positive root
of the equation
$$
2h^2\exp(1+\sqrt{1+4h^2}-2h)= 1 + \sqrt{1+4h^2}
$$
and plays a key role in the following result (which is proved in
Section 2):
\begin{lem} \label{c1} Let $\epsilon \in (0,0.25], \ h >0$. Then
the characteristic function $\psi(z,\epsilon):=\epsilon
z^2-z-\exp(-z h)$ has exactly two (counting multiplicity) negative
zeros $\lambda_1\leq \lambda_2<0$ if and only if one of the
following conditions holds
\begin{enumerate}
  \item $0< h \leq 1/e $,
  \item $ \epsilon \geq \epsilon^*(h)$ and $1/e < h \leq  \ h_1$.
\end{enumerate}
\end{lem}
\noindent Here the continuous $\epsilon^*(h)$ is defined in
parametric form by
$$
\epsilon^*(h(t))= th(t), \ h(t) =
(2t+\sqrt{4t^2+1})\exp(-1-\frac{2t}{1+\sqrt{4t^2+1}}), \ t \in [0,
0.445\dots].
$$
Let us state now the main results of this paper.
\begin{thm} \label{main} Eq. (\ref{17}) has a positive monotone  wavefront
$u = \varphi(\nu\cdot x + ct), $ \ $ |\nu|=1,$ connecting $0$ with
$1$ if and only if one of the following conditions holds
\begin{enumerate}
  \item $0\leq h \leq 1/e = 0.367879441...$ and $2 \leq c < c^*(h):=+\infty$;
  \item $1/e < h \leq  \ h_1 =
  0.560771160\dots$ and $2 \leq c \leq c^*(h):=
  1/\sqrt{\epsilon^*(h)}$.
\end{enumerate}
Furthermore, set $\phi(s):=\varphi(cs)$.  Then for some appropriate
$\phi_-$ (given below explicitly), we have that $\phi = \lim_{j \to
+\infty} \mathcal{A}^j\phi_-$ (if $c>2$), and $\phi = \lim_{j \to
+\infty} \mathcal{B}^j\phi_-$ (if $c=2$), where the convergence is
monotone and uniform on $\R$. Finally, for each fixed
$c\not=c^*(h)$, $\phi(t)$ is the only possible profile (modulo
translation) and $\phi(t), \phi_-(t)$ have the same asymptotic
representation $1-e^{\lambda_2t}(1+o(1))$ at $+\infty$.
\end{thm}
\vspace{-1mm}
\begin{cor}
If $h > h_1=0.560771160\dots$ then the delayed KPP-Fisher equation
does not have any positive monotone traveling wavefront.
\end{cor}
\vspace{-2mm} Next, let us define the continuous function
$\epsilon^\#(h)$ parametrically by
\begin{equation}\label{eg}\hspace{-10mm}
\epsilon^\#(h(t))= \frac{t+2+\sqrt{2t+4}}{t^2},  h(t) =
-\frac{\ln(2+\sqrt{2t+4})}{t}, \ t \in (-2,-1.806\dots]
\end{equation}
\vspace{-4mm} Set $h_0: = 0.5336619208 \dots$ (see also Lema
\ref{c1dd} for its complete definition) and
$$
c^\#(h):= \left\{%
\begin{array}{lll}
 +\infty ,
    & \hbox{when}\  h \in(0,0.5\ln2],
   \\
1/\sqrt{\epsilon^\#(h)},
    & \hbox{when} \ h \in(0.5\ln2, h_0],   \\
      2, &  \hbox{when}\
      h > h_0.
\end{array}%
\right.\nonumber $$\vspace{-2mm}
\begin{figure}[h]
\centering \fbox{\includegraphics[width=6cm]{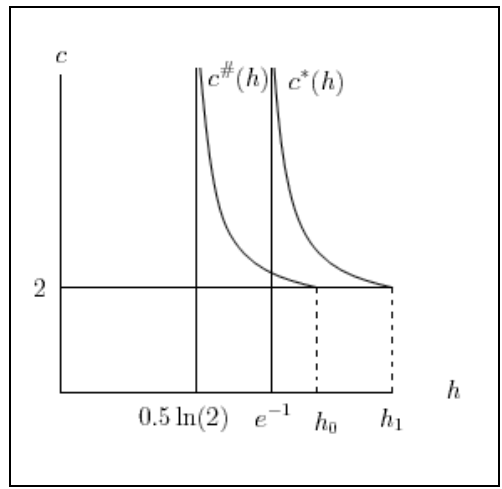}}
\caption{Schematic presentation of the critical speeds  and
delays.}\end{figure} \vspace{3mm}

\begin{thm} \label{mainas} Let $u = \varphi(\nu\cdot x + ct), $ \ $ |\nu|=1,$ be a positive monotone
traveling front of Eq. (\ref{17}).  Set $\phi(s):=\varphi(cs)$.
Then, for some appropriate $t_0$, positive $K_j$ and every small
positive $\sigma$, we have at $t = -\infty$\newpage
$$\hspace{-10mm}\nonumber
\phi(t+t_0)=\left\{%
\begin{array}{lll}\nonumber -K_2te^{\lambda t} + O(e^{(2\lambda-\sigma) t}), &  \hbox{when}\
      c=2,
    \\ \nonumber
     e^{\lambda t} - K_1 e^{\mu t} + O(e^{(2\lambda-\sigma) t}),
    & \hbox{when}\, \ 2<c< 1.5\sqrt{2}, \\ \nonumber
     e^{\lambda t}  + O(e^{(2\lambda-\sigma) t}),
    & \hbox{when}\  c\geq 1.5\sqrt{2}=2.121\dots
\end{array}%
\right.\nonumber $$ Similarly, at $t = +\infty$
\begin{equation}\label{3up}\hspace{-10mm}
\phi(t+t_0)=\left\{%
\begin{array}{llll} \label{asfor}
       1-e^{\lambda_2 t}  + O(e^{(2\lambda_2+\sigma) t}),
    & \hbox{when \ }\ h \leq h_0, \ c\in [2, c^\#(h)]\cap \R,    \\
     1-e^{\lambda_2 t} +K_3e^{\lambda_1 t}+
    & \hbox{when }\ h \in (0.5\ln2, h_1] \     \\
   \hspace{5mm}+O(e^{(\lambda_1-\sigma) t}), & {\hbox{and }  c \in (c^\#(h), c^*(h))},   \\
      1-K_4te^{\lambda_2 t} + O(e^{(\lambda_2-\sigma) t}), &  \hbox{when}\
      c=c^*(h)\ \hbox{and } h \in (1/e,h_1].
\end{array}%
\right.\nonumber
\end{equation}
\end{thm}
Theorem \ref{mainas} suggests the way of  approximating the
traveling front profile: e.g., for $c\not=2,c^*(h)$, we can take
functions $a_-(t):= c_1e^{-\lambda t}$ and $a_+(t):= 1- e^{\lambda_2
t}$ and glue them together at some point $\tau$.  The point $\tau$
and $c_1>0$ have to be chosen to assure maximal smoothness of the
approximation at $\tau$. As we will see in Section 3, this idea
allows to construct reasonable lower approximations to the exact
traveling wave. See also Figure 2 below.

\begin{rem} As it was showed by Ablowitz and Zeppetella
\cite{AZ},  equation (\ref{kppf}) has the explicit exact wavefront
solution $u = \varphi_{\star}(\nu\cdot x + ct), $ \ $ |\nu|=1,$ with
$c=5/\sqrt{6}=2.041\dots$ and the (scaled) profile
$$
\phi_{\star}(s)=
\left(\frac{1}{2}+\frac{1}{2}\tanh(\frac{5s}{12}+s_0) \right)^2, \
\phi_{\star}(s):=\varphi_{\star}(cs).
$$
If we select $s_0=0.5\ln2$, then $$\phi_{\star}(s)=1 -
2e^{-5s/6-2s_0}+O(e^{-5s/4})= 1 - e^{-5s/6}+O(e^{-5s/4}),\ s \to
+\infty,$$ so that  $\phi_{\star} = \lim_{j \to +\infty}
\mathcal{A}^j\phi_-$ in view of Theorem \ref{main} and the
uniqueness (up to translations) of the traveling front for the
non-delayed KPP-Fisher equation. Figure 2 (on the left) shows five
approximations $\mathcal{A}^j\phi_-, j =0,1,2,3,4,$ and the exact
solution $\phi_{\star}$, the graphs are ordered as $\phi_- <
\mathcal{A}\phi_-< \mathcal{A}^2\phi_-< \mathcal{A}^3\phi_-<
\phi_{\star}$. On the right, the four first approximations
$\mathcal{B}^j\phi_-, j =0,1,2,3,$ of $\phi$ are plotted when $c=2,
\ h=0.56$. It should be noted that the limit function $\phi$ and the
initial approximation $\phi_-$ have the same first two terms
$(1-\exp(\lambda_2 t))$ of their asymptotic expansions at $+\infty$.
See Theorem \ref{main} and Sections 3,4. However, as the analysis of
the Ablowitz-Zeppetella solution shows, these $\phi$ and  $\phi_-$
may have different first terms of their expansions at $-\infty$.
This partially explains a better agreement between the exact
solution and their approximations for $t \geq \tau= 0.487\dots$ on
the left picture (the value of $\tau$ is given in Section 3).
\end{rem}
\begin{figure}[h]\label{ffig}
\centering \fbox{\includegraphics[width=6.3cm]{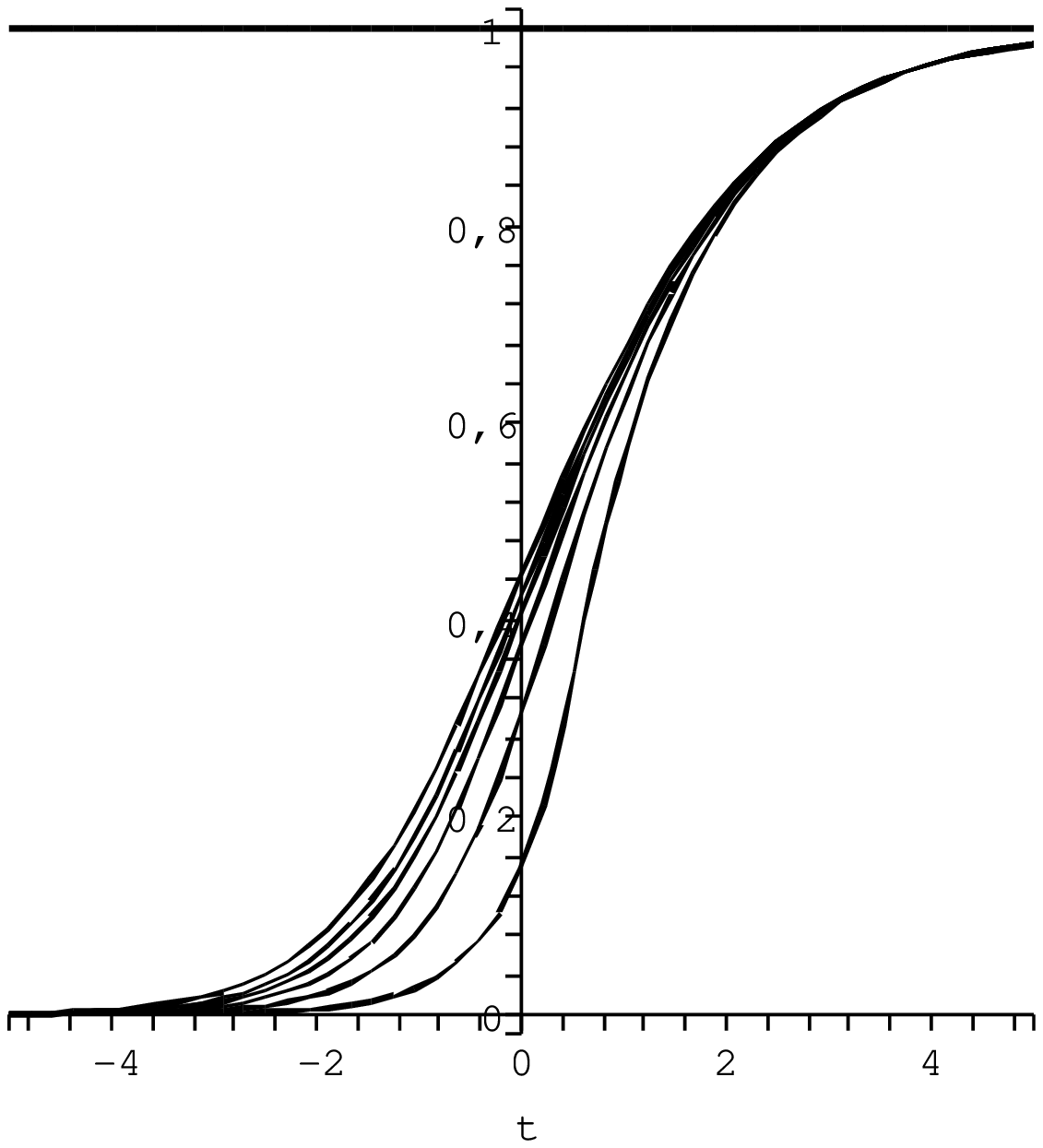}}
\fbox{\includegraphics[width=6.3cm]{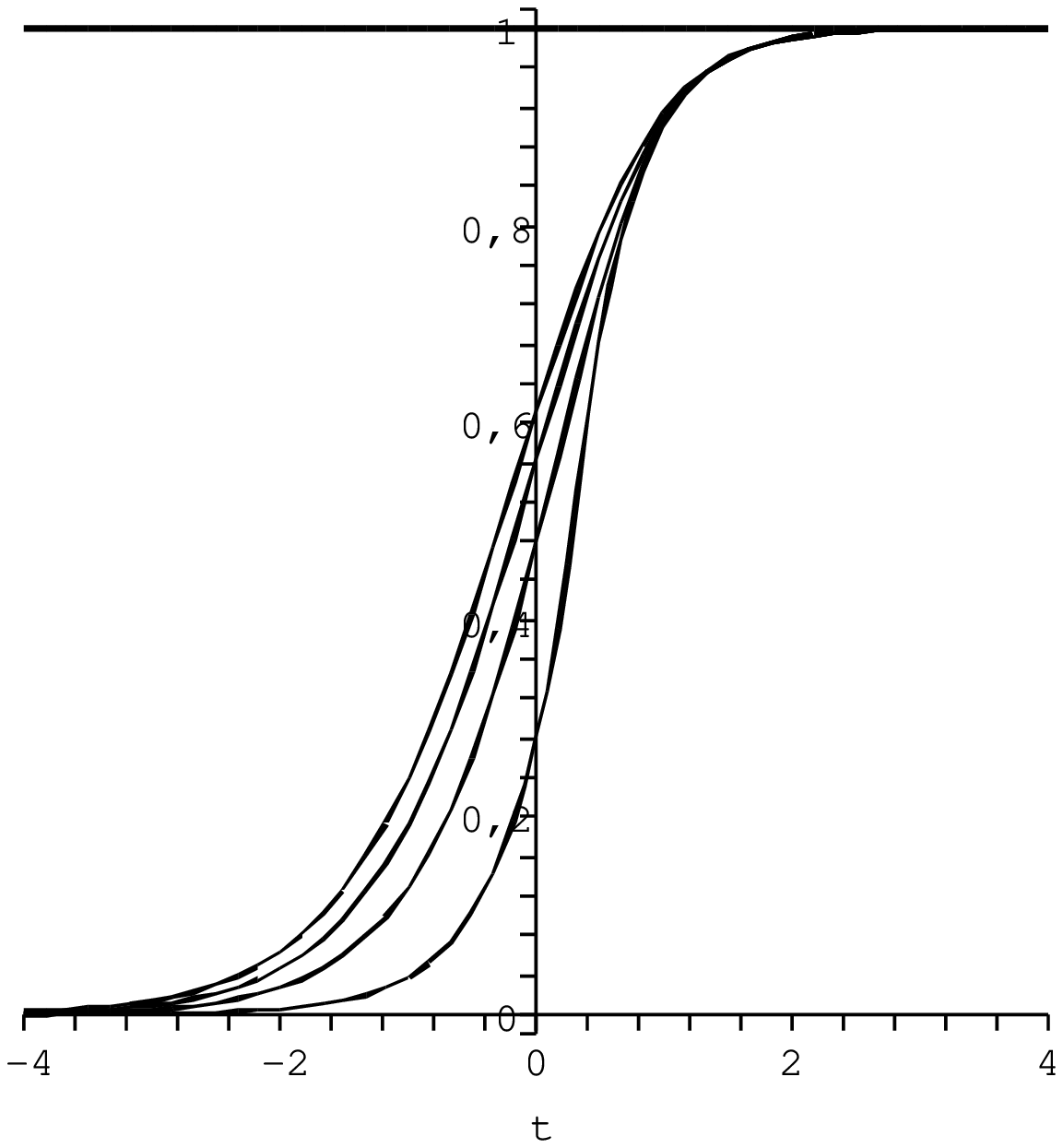}} \caption{On the
left: increasing sequence of approximated waves
$\mathcal{A}^j\phi_-, j =0,1,2,3,4,$ and the Ablowitz-Zeppetella
exact solution $\phi_{\star}$ ($\epsilon= 0.24$ and $h=0$). On the
right: approximations $\mathcal{B}^j\phi_-, j =0,1,2,3$ ($\epsilon=
0.25$ and $h=0.56$). }\end{figure} The structure of the remainder of
this paper is as follows. In Section 2, the characteristic function
of the variational equation at the positive steady state is
analyzed. In the third [the fourth] section, we present a lower [an
upper] solution. Section \ref{comments} contains some comments on
the smoothness of upper and lower solutions. Theorems \ref{main} and
\ref{mainas} are proved in Sections  \ref{sm} and \ref{rrra},
respectively.
\section{Characteristic equation at the positive steady state}
In this section, we study the zeros of $\psi(z,\epsilon):=\epsilon
z^2-z-\exp(-z h), \ \epsilon,h >0.$ It is straightforward  to see
that $\psi$ always has a unique positive simple zero. Since
$\psi'''(z,\epsilon)$ is positive, $\psi$ can have at most three
(counting multiplicities) real zeros, one of them positive and the
other two (when they exist) negative. Lemma \ref{c1} in the
introduction provides a criterion for the existence of two negative
zeros $\lambda_1\leq \lambda_2<0$. We start by proving this
result:\vspace{0mm}
\begin{pf}[Lemma \ref{c1}]
Consider the equation $-z=\exp(-z h).$ An easy  analysis shows that
(i) this equation has exactly two real simple solutions $z_1<z_2<0,$
$ z_2 > -e,$ if $h \in (0,1/e)$, (ii) it has one double real root
$z_1=z_2=-e$ if $h=1/e$, and (iii) it does not have any real root if
$h> 1/e$. As a consequence,
\begin{equation}\label{g2}
\epsilon z^2-z=\exp(-z h)
\end{equation}
\vspace{-5mm} has two negative simple solutions if $\epsilon >0$ and
$h \in (0,1/e]$.

A similar argument shows that for every $h > 1/e$ there exists
$\epsilon^*(h)>0$ such that Eq. (\ref{g2}) (a) has two negative
simple roots if $\epsilon
>\epsilon^*(h)$, (b) has one negative double root if $\epsilon
=\epsilon^*(h)$, (c) does not have any solution if $\epsilon
<\epsilon^*(h)$. In particular, $\epsilon=\epsilon^*(h), z=
\lambda_1(h)=\lambda_2(h),$ solve the system
$$
\epsilon z^2-z=\exp(-z h), \quad 2\epsilon z-1=-h\exp(-z h),
$$
which yields the parametric representation for $\epsilon^*(h)$ given
in the introduction.

Finally, a direct graphical analysis of (\ref{g2}) shows that
$\epsilon^*(h)$ is increasing with respect to  $h$. Hence, since
$\epsilon^*(h)\leq 0.25$, we conclude that $h \leq
(\epsilon^*)^{-1}(0.25)=: h_1= 0.560771\dots$ \hfill
{$\square$}\end{pf}
\begin{lem}\label{c1dd} Let $\lambda_1\leq
\lambda_2<0$ be two negative zeros of $\psi(z,\epsilon)$ and
$\epsilon \in (0,0.25]$ be fixed. Then  $\lambda_1 \leq 2\lambda_2$
if and only if one of the following conditions holds
\begin{enumerate}
  \item $0 <h \leq 0.5\ln 2 = 0.347\dots$;
  \item $ \epsilon \geq \epsilon^\#(h)$ and $ 0.5\ln 2 < h \leq  \ h_0:= 0.5336619208\dots$.
\end{enumerate}
\end{lem}\begin{pf} This lemma can be proved analogously to the previous
one, we briefly outline the main arguments. First, for each fixed
positive $\epsilon^\#$ we may find $h(\epsilon^\#)>0$ such that
$\lambda_1 < 2\lambda_2$ if $h \in (0,h(\epsilon^\#))$ and
$\lambda_1 = 2\lambda_2$ if $h=h(\epsilon^\#)$. In this way,
$$
\epsilon^\# \lambda_2^2-\lambda_2=\exp(-\lambda_2 h(\epsilon^\#)),
\quad 4\epsilon^\# \lambda_2^2-2\lambda_2=\exp(-2\lambda_2
h(\epsilon^\#)),
$$
which yields representation (\ref{eg}). Now, we complete the proof
by noting that $h(\epsilon)$ is continuous  and strictly increasing
on $(0,+\infty)$ and $h(0+)= 0.5\ln2, \ h_0 = h(0.25)$. \hfill
{$\square$}\end{pf}
\begin{lem}\label{c1de} Let $\lambda_1\leq
\lambda_2<0$ be two negative zeros of $\psi(z,\epsilon)$ and
$\epsilon \in (0,0.25]$ be fixed. Then  $\Re\lambda_j < \lambda_1$
for every complex root of $\psi(z,\epsilon)=0$.
\end{lem}
\vspace{-5mm}
\begin{pf}
Set $\alpha: = (1+2\epsilon -\sqrt{1+4\epsilon^2})/(2\epsilon),
 \ a:= -e^{-\alpha
h}/(\sqrt{1+4\epsilon^2}-2\epsilon), \
k:=\epsilon/(\sqrt{1+4\epsilon^2}-2\epsilon)$. Then $\alpha, k
>0, a <0,$ and
$$\psi(z+\alpha) = (\sqrt{1+4\epsilon^2}-2\epsilon)(kz^2-z-1+ae^{-zh}).$$
It is easy to see that $p(z):=kz^2-z-1+ae^{-zh}$ also has two
negative and one positive root. Since the translation $z \to
z+\alpha$ of the complex plain does not change the mutual position
of zeros of $\psi$, the statement of Lemma \ref{c1de} follows now
from \cite[Remarks 19,20]{TT}. \hfill {$\square$}\end{pf}
\vspace{-9mm}

\section{A lower solution when $\lambda_1 < \lambda_2$}\label{low1}
\vspace{-3mm} In this section, we assume either condition (1) or
condition (2) of Theorem \ref{main} holds. In addition, let $c \in
[2,c^*(h))$ so that $\lambda_1 < \lambda_2$ (where
$\lambda_1:=-\infty$ if $h=0$) and $\lambda\leq \mu$. Set
\vspace{-5mm}
$$\tau = \frac{1}{\lambda_2}\ln \frac{\lambda}{\lambda- \lambda_2}
>0,\quad
\phi_-(t)=\left\{%
\begin{array}{ll}{
     \frac{-\lambda_2}{\lambda- \lambda_2}e^{\lambda(t- \tau)}},
    & \hbox{if}\ t \leq \tau,   \vspace{5mm} \\
     1 - e^{\lambda_2 t} & \hbox{if} \
    t \geq \tau .
\end{array}%
\right.
$$
It is easy to see that $\phi_- \in C^1(\mathbb{R})\cap
C^2(\mathbb{R}\setminus\{\tau\})$ with $\phi_-'(t) >0, \ t \in
\mathbb{R}$, and
\begin{equation} \label{ls}
\epsilon \phi_-''(t) - \phi_-'(t) + \phi_-(t)(1- \phi_-(t-h)) < 0,
\quad t \in \R \setminus (\tau, \tau +h].
\end{equation}
\begin{lem}\label{low} Inequality (\ref{ls}) holds for all $t \in \R$.
\end{lem}\vspace{-5mm}\begin{pf} The case $h=0$ is obvious, so let $h >0$. It suffices to consider $t \in (\tau, \tau
+h]$. If we take $t \in (\tau, \tau +h]$, then
$$
\epsilon \phi_-''(t) - \phi_-'(t) + \phi_-(t)(1- \phi_-(t-h)) =
-\epsilon \lambda_2^2e^{\lambda_2 t} + \lambda_2e^{\lambda_2 t} + $$
$$
( 1 - e^{\lambda_2 t})(1+ \frac{\lambda_2}{\lambda-
\lambda_2}e^{\lambda(t- \tau-h)}) = - e^{\lambda_2 (t-h)} + ( 1 -
e^{\lambda_2 t})(1+ \frac{\lambda_2}{\lambda-
\lambda_2}e^{\lambda(t- \tau-h)}) =
$$
$$
1 - e^{\lambda_2 (t-h)} + \frac{\lambda_2}{\lambda-
\lambda_2}e^{\lambda(t- \tau-h)}- e^{\lambda_2 t}  - e^{\lambda_2 t}
\frac{\lambda_2}{\lambda- \lambda_2}e^{\lambda(t- \tau-h)} =
$$
$$
1 + \frac{\lambda_2}{\lambda- \lambda_2}e^{\lambda s} -
\frac{\lambda}{\lambda- \lambda_2} e^{\lambda_2 s} -
\frac{\lambda}{\lambda- \lambda_2}e^{\lambda_2 (s + h)} -
\frac{\lambda\lambda_2}{(\lambda- \lambda_2)^2} e^{\lambda_2 (s +
h)} e^{\lambda s} =: \rho(s)
$$
where $ s = t- \tau-h \in (-h,0]$. The direct differentiation shows
that
$$
\rho'(s) =  \frac{-\lambda_2\lambda}{\lambda-
\lambda_2}\left[-e^{\lambda s} +e^{\lambda_2 s} + e^{\lambda_2 (s +
h)}(1 + \frac{\lambda + \lambda_2}{\lambda- \lambda_2} e^{\lambda
s})\right] > 0,
$$
$$ \hspace{-25mm}{\rm since}\  e^{\lambda s} \leq 1, \ e^{\lambda_2 s}\geq 1,\ {\rm
and}\ (1 + \frac{\lambda+ \lambda_2}{\lambda- \lambda_2} e^{\lambda
s})
> 1, {\rm if } \ \lambda+ \lambda_2 \geq 0, $$
$$
(1 + \frac{\lambda+ \lambda_2}{\lambda- \lambda_2} e^{\lambda s})
\geq 1 + \frac{\lambda+ \lambda_2}{\lambda- \lambda_2} =
\frac{2\lambda}{\lambda- \lambda_2} >0 , {\rm if } \ \lambda+
\lambda_2 <  0. $$ Finally, we have that $\rho(s) < 0$ for all $s
\in [-h,0]$ since $\rho'(s) > 0$ and
$$
\rho(0)=   - \lambda^2(\lambda- \lambda_2)^{-2} e^{\lambda_2 h}
<0.\hspace{5cm} \square
$$
\end{pf}
\begin{rem}[A lower solution when $\lambda_1 = \lambda_2$]
We can not use  $\phi_-$ as a lower solution  when $c = c^*(h), \
1/e < h \leq   \ h_1$. Indeed, by Theorem \ref{mainas}, in this case
$\phi_-$ converges to the positive steady state faster than the
heteroclinic solutions. In Section \ref{comments}, we will present
an adequate lower solution for this situation. However, it will not
be $C^1$-smooth.
\end{rem}
\vspace{-5mm}
\section{An upper solution when $\lambda_1 < \lambda_2$} \label{upp}
\vspace{-5mm} Suppose that $\lambda_1 < \lambda_2$ and  set
$\phi_2(t): = 1 - e^{\lambda_2 t} + e^{ r t}$ for some $r \in
(\lambda_1, \lambda_2)$. Recall that $\lambda_1:=-\infty$ if $h=0$.
Obviously, $\psi(r,\epsilon)
>0$ and $\phi_2(t) \in (0,1)$ for $t
> 0$. Next, it is immediate to check that
$\phi_2:\R\to\R$ has a unique critical point (absolute minimum) $t_0
= t_0(r)>0$:
$$
t_0(r)= \frac{\ln (-r) - \ln (-\lambda_2)}{\lambda_2-r}, \quad
\lambda_2
 e^{\lambda_2 t_0} = r e^{ r t_0}.
$$
Observe that if $h \in (0, 1/e)$, then we can assume that $t_0(r)
\geq h$ since
$$
\lim_{r \to \lambda_2-}t_0(r) = - 1/\lambda_2 > 1/e> h,
$$
where the last inequalities were established in the proof of Lemma
\ref{c1}. It is clear that the function
$$
\phi_+(t)=\left\{%
\begin{array}{ll}
     \phi_2(t),
    & \hbox{if}\ t \geq t_0(r),   \\
     \phi_2(t_0(r)), & \hbox{if} \ t \leq t_0(r)
\end{array}%
\right.
$$
is $C^1$-continuous and increasing on $\R$. Moreover, $\phi_+(t)\in
C^2(\R\setminus \{t_0(r)\})$.
\begin{lem}\label{up} For all $r < \lambda_2$ sufficiently close to
$\lambda_2$, $\phi_+$ satisfies the inequality
$$
\epsilon \phi''(t) - \phi'(t) + \phi(t)(1- \phi(t-h)) \geq 0, \quad
t \in \R. $$
\end{lem} \vspace{-5mm}
\begin{pf} Step I. First we prove that,  for all $t
\geq t_0$, the following inequality holds:
$$
(\mathfrak{N}\phi_2)(t):=\epsilon \phi''_2(t) - \phi'_2(t) +
\phi_2(t)(1- \phi_2(t-h)) \geq 0.
$$
In particular, this implies that $(\mathfrak{N}\phi_+)(t)\geq 0$ if
$t \geq t_0+h$.  For $t = t_0 +s$, we have that
$$
(\mathfrak{N}\phi_2)(t) = \psi(r,\epsilon)e^{ r
t}-\psi(\lambda_2,\epsilon)e^{\lambda_2 t} + (-e^{\lambda_2 t} + e^{
r t})(e^{\lambda_2 (t-h)} - e^{ r (t-h)})=
$$
$$
\psi(r,\epsilon)e^{ r t} + (-e^{\lambda_2 t} + e^{ r
t})(e^{\lambda_2 (t-h)} - e^{ r (t-h)})=
$$
$$
e^{ r t_0}\left[\psi(r,\epsilon)e^{ r s} +
(-\frac{r}{\lambda_2}e^{\lambda_2 s} + e^{ r s})e^{ r
t_0}(\frac{r}{\lambda_2}e^{\lambda_2 (s-h)} - e^{ r (s-h)})\right]=
$$
$$
e^{ r (t_0+s)}\left[\psi(r,\epsilon) +
(-\frac{r}{\lambda_2}e^{(\lambda_2-0.5r) s} + e^{ 0.5r s})e^{ r
t_0}(\frac{r}{\lambda_2}e^{-\lambda_2 h}e^{(\lambda_2 -0.5r)s} - e^{
-rh}e^{ 0.5r s})\right]=
$$
$$
e^{ rt}\left[\psi(r,\epsilon) + A_1(s)e^{ r t_0}A_2(s)\right].
$$
It is easy to see that $A_j(+\infty) =0$ and that $A_j$ has a unique
critical point $s_j$, with
$$
\lim_{r \to \lambda_2-} s_1(r) = -1/\lambda_2, \lim_{r \to
\lambda_2-} s_2(r) = h-1/\lambda_2.
$$
Therefore, for some small $\delta >0$ and for all $r$ close to
$\lambda_2$, the function $A_1(s)e^{ r t_0}A_2(s)$ is strictly
increasing to $0$ on the interval $[ h-1/\lambda_2+\delta, +\infty)$
and it is strictly decreasing on $[0, -1/\lambda_2-\delta]$. This
means that if $(\mathfrak{N}\phi_2)(t) \geq 0$ for all $t \in
[t_0-1/\lambda_2-\delta, t_0 + h-1/\lambda_2+\delta]$ then
$(\mathfrak{N}\phi_2)(t) \geq 0$ for $ t \geq t_0$. In order to
prove the former, consider the expression
$$
\frac{e^{-r t_0}}{r-\lambda_2}\left(\epsilon \phi''_2(t) -
\phi'_2(t) + \phi_2(t)(1- \phi_2(t-h))\right)= $$ $$
\frac{\psi(r,\epsilon)e^{ r s} + (-\frac{r}{\lambda_2}e^{\lambda_2
s} + e^{ r s})e^{ r t_0}(\frac{r}{\lambda_2}e^{\lambda_2 (s-h)} -
e^{ r (s-h)})}{r-\lambda_2}: = \Gamma_{\epsilon}(r,s).
$$
Since $\Gamma_{\epsilon}(r,s)$ is analytical on some open
neighborhood $\Omega \subset \R^2$ of the compact segment
$\{\lambda_2\}\times [-1/\lambda_2-\delta, h-1/\lambda_2+\delta]
\subset \R^2$, we find that, for every fixed $\epsilon >0$,
$$ \lim_{r\to \lambda_2-}\Gamma_{\epsilon}(r,s)
= \psi'(\lambda_2,\epsilon)e^{ \lambda_2 s} < 0
$$
uniformly on $[-1/\lambda_2-\delta, h-1/\lambda_2+\delta]$. As a
consequence, we obtain that
$$
\epsilon \phi''_2(t) - \phi'_2(t) + \phi_2(t)(1- \phi_2(t-h)) > 0,
\quad  t \in [t_0-1/\lambda_2-\delta, t_0+ h-1/\lambda_2+\delta].
$$
Step II. Now, we are ready to prove that $(\mathfrak{N}\phi_+)(t)
\geq 0, \ t \in [t_0, t_0+h].$ Indeed, since $\phi_2(t_0) \leq
\phi_2(t-h)$ for $t \in [t_0, t_0+h]$, we have that
$$(\mathfrak{N}\phi_+)(t)= \epsilon \phi_2''(t) - \phi_2'(t) + \phi_2(t)(1- \phi_2(t_0)) \geq
$$
$$
\phi_2''(t) - \phi_2'(t) + \phi_2(t)(1- \phi_2(t-h)) \geq 0, \quad t
\in [t_0, t_0+h].
$$
Finally, since the inequality $(\mathfrak{N}\phi_+)(t) > 0, \ t \leq
t_0$, is obvious, the proof of the lemma is completed. \hfill
{$\square$}\end{pf}

\begin{rem}[An upper solution when $\lambda_1 = \lambda_2$]
We can not use  $\phi_+$ as an upper solution  when $c = c^*(h), \
1/e < h \leq   \ h_1$. Moreover, in this case it is not difficult to
show that $\phi_+$  satisfies inequality (\ref{ls}) for all $r <
\lambda_2$ sufficiently close to $\lambda_2$ and for large positive
$t$.
\end{rem}
\vspace{-5mm}
\section{Some comments on upper and lower solutions}
\vspace{-5mm} \label{comments}
\subsection{Non-smooth solutions}
\vspace{-5mm} The problem of smoothness of the lower (upper)
solutions is an interesting and important aspect of the topic, see
\cite{bn,ma}. As we have seen in the previous sections,
$C^1-$smoothness condition can be rather restrictive even when a
simple nonlinearity (the birth function) is considered. The above
mentioned works \cite{ma} show that continuous and piece-wise
$C^1-$continuous lower (upper) solutions $\phi_{\pm}$ still can be
used if some sign conditions are fulfilled at the points of
discontinuity of $\phi_{\pm}'$. Moreover, as we prove it below even
discontinuous functions $\phi_{\pm}$ can be also used. We start with
a simple result of the theory of impulsive systems \cite{sp} which
can be viewed as a version of the Perron theorem for piece-wise
continuous solutions, cf. \cite{bn}.
\begin{lem}\label{imp} Let $\psi: \R\to \R$ be a bounded classical solution of
the second order  impulsive equation
$$
\psi'' + a\psi' + b\psi =f(t), \quad \Delta \psi|_{t_j} = \alpha_j,
\quad \Delta \psi'|_{t_j} = \beta_j,
$$
where $\{t_j\}$ is a finite increasing sequence, $f: \R\to \R$ is
bounded and continuous at every $t \not= t_j$ and the operator
$\Delta$ is defined by $\Delta w|_{t_j} := w(t_j+)-w(t_j-)$. Assume
that $z^2 + az + b =0$ has two positive roots $0<\lambda\leq\mu$.
Then
\begin{eqnarray}\label{pfo}&&\hspace{-10mm}
 {\rm \ if } \ \lambda< \mu \ {\rm\ we \ have \ that }\ \psi(t)=\frac{1}{\mu - \lambda} \int_t^{+\infty}\left(e^{\lambda
(t-s)}- e^{\mu (t-s)}\right)f(s)ds\\ &&  \nonumber
 +\frac{1}{\mu - \lambda}\sum_{t<t_j}\left[\left(\lambda e^{\mu(t-t_j)}-\mu e^{\lambda(t-t_j)}\right)\alpha_j+
 \left(e^{\lambda(t-t_j)}-e^{\mu(t-t_j)}\right)\beta_j\right], \
 \ t\not=t_j;\\ && \hspace{-10mm} \nonumber  {\rm \ if } \ \lambda= \mu=-0.5a \ {\rm\ we \ have \ that}\
 \psi(t) =  \int_t^{+\infty}(s-t)e^{-0.5a(t-s)}f(s)ds+\\  \nonumber &&
\sum_{t<t_j}e^{-0.5a(t-t_j)} \left[(t_j-t)(\beta_j+0.5a\alpha_j) -
\alpha_j\right], \  \ t\not=t_j.
\end{eqnarray}
\end{lem}\vspace{-10mm}\begin{pf}
See \cite[Theorem 87]{sp}. \hfill {$\square$}\end{pf} \vspace{-5mm}
Next, the corollary below shows that our lower solution is an upper
solution in the sense of Wu and Zou \cite{wz}:
\begin{cor} \label{cls} Assume that $\psi: \R\to \R$ is bounded and such that
the derivatives $\psi',\psi'': \R\setminus\{t_j\}\to \R$ exist and
are bounded. Suppose also that $\psi$ is a classical solution of the
impulsive inequality $$ \psi'' + a\psi' + b\psi \leq f(t), \quad
\Delta \psi|_{t_j} = \alpha_j, \quad \Delta \psi'|_{t_j} = \beta_j.
$$
If  $\alpha_j \geq 0, \ \beta_j \leq 0$,  then \vspace{-5mm}
\begin{eqnarray}
\nonumber &&   \psi(t) \leq \frac{1}{\mu - \lambda}
\int_t^{+\infty}\left(e^{\lambda (t-s)}- e^{\mu (t-s)}\right)f(s)ds,
\quad {\rm when}\ \lambda< \mu, \\
\nonumber &&  \psi(t) \leq
\int_t^{+\infty}(s-t)e^{-0.5a(t-s)}f(s)ds, \quad {\rm when}\
\lambda= \mu=-0.5a.
\end{eqnarray}
\vspace{-5mm}
\end{cor}
\vspace{-5mm}
\begin{pf} Suppose that $\lambda<\mu$, the case  $\lambda=\mu$ is similar. Clearly, $q(t):= f(t)- (\psi''(t) + a\psi'(t) + b\psi(t)) \geq 0$
and $\lambda e^{\mu(t-t_j)}<\mu e^{\lambda(t-t_j)},$ $\
e^{\lambda(t-t_j)}> e^{\mu(t-t_j)}$ for $t < t_j$. Thus the desired
inequality follows  from (\ref{pfo}). \hfill {$\square$}\end{pf}
\vspace{-7mm}

\subsection{A lower solution when $\lambda_1 = \lambda_2, \ h \in (1/e,
h_1]$}\vspace{-5mm}\label{lsc} In Section \ref{low1}, a lower $C^1-$
solution was presented for the case when $\lambda_1 < \lambda_2$.
However, to apply our iterative procedure in the critical case
$\lambda_1 = \lambda_2$, we also need to construct a lower solution
for the corresponding range of parameters. It is worth to mention
that our approach does not require any upper solution once a lower
solution is found and the existence of the heteroclinic is proved,
see Corollary \ref{mcoro}. Here, we provide a continuous and
piece-wise analytic lower solution $\phi_-(t)$ if $\lambda_1 =
\lambda_2$. Our solution has a unique singular point $\tau'$ where
$\Delta \phi_-|_{\tau'} = 0, \quad \Delta \phi_-'|_{\tau'}
>0$. This shows that, in general, the sign conditions of Corollary \ref{cls}
need not to be satisfied.

Take some positive $ A > (e^{-\lambda_2 h}-1)/h$ and let $\tau'$ be
the positive root of the equation $ At + 1 = e^{-\lambda_2 t}$. It
is easy to see that $\tau' >h$. Consider the piece-wise smooth
function $\phi_-: \R \to [0,1)$ defined by
\begin{equation}\label{sokr}
\phi_-(t)=\left\{%
\begin{array}{ll}
    0, & \hbox{if} \ t \leq \tau',    \\
     1 - (At+1)e^{\lambda_2 t}, & \hbox{if} \
    t \geq \tau'.
\end{array}%
\right.
\end{equation}
\begin{prop}\label{mcr} The inequality $(\mathcal{K}\phi_-)(t) >
\phi_-(t)$ holds for all $t \in \R$.
\end{prop}\vspace{-5mm}\begin{pf} Below, we are assuming that $h\not= h_1$ so that $\lambda <
\mu$ and $\mathcal{K}= \mathcal{A}$; however, a similar argument
works also in the case  $h=h_1$ (when $\mathcal{K}= \mathcal{B}$).
It suffices to prove that $(\mathcal{A}\phi_-)(t) > \phi_-(t)$ for
$t \geq \tau'$. Let $C^2-$ smooth function $\psi$ be defined  by
$$
\psi(t)=\left\{%
\begin{array}{lll}
    1-(At+1)e^{\lambda_2 t}, & \hbox{if} \ t \geq \tau'-h ,   \\
     B(t), & \hbox{if} \
    0 \leq t \leq \tau'-h ,    \\
    0, & \hbox{if} \ t \leq 0,
\end{array}%
\right.
$$
for some appropriate  continuous decreasing  $B(t)$. Set
$$
\zeta(t): = \epsilon \psi''(t) - \psi'(t) + \psi(t)(1- \psi(t-h)).
$$
It is easy to check that $\zeta \in C(\R,\R)$ is bounded on $\R$ and
$\zeta(t) < 0$ for all $t > \tau'$. But then, for all $t
> \tau'$, we have that
$$
\phi_-(t) = \psi(t) = (\mathcal{A}\psi)(t)+ \frac{1}{\epsilon(\mu -
\lambda)} \int_t^{+\infty}\left(e^{\lambda (t-s)}- e^{\mu
(t-s)}\right)\zeta(s)ds < (\mathcal{A}\psi)(t)\leq
$$
$$
= \frac{1}{\epsilon (\mu - \lambda)} \int_t^{+\infty}(e^{\lambda
(t-s)}- e^{\mu (t-s)})\phi_-(s)\phi_-(s-h)ds =
(\mathcal{A}\phi_-)(t).\quad \square
$$
\end{pf}\vspace{-5mm}
\subsection{Ordering the upper and lower solutions}\vspace{-5mm} Finally, we
show that the condition of the correct ordering $\phi_- \leq \phi_+$
is not at all restrictive provided that solutions $\phi_{\pm}$ are
monotone and satisfy some natural asymptotic relations.

\begin{lem}\label{as} Assume that functions $\phi_{\pm}:\R \to [0,1),$ $ j =1,2,$ are increasing
and, for some fixed $k \in \{0,1\}$, the following holds
$$
 \lim_{t \to
-\infty}\phi_{\pm}(t)e^{-\lambda t} = \alpha_\pm,\ \lim_{t \to
+\infty}(1-\phi_{\pm}(t))t^{-k}e^{-\lambda_2 t} = \beta_{\pm,k},
$$
where $\beta_{\pm,k}>0$,  and $\alpha_{-}\in [0,+\infty)$,
$\alpha_{+}\in (0,+\infty]$. Then there exists a real number
$\sigma$ such that $\phi_-(t) < \phi_+(t+ \sigma)$ for all $t \in
\R$.
\end{lem}\vspace{-5mm} \begin{pf} It is clear that  $\phi_{-}(-\infty)=0$ and  $\phi_{\pm}(+\infty) =1$.
Let $\sigma_0$ be sufficiently large to satisfy
$\beta_{+,k}e^{\lambda_2 \sigma_0} < \beta_{-,k}, \ \alpha_-
e^{-\lambda \sigma_0} <\alpha_{+}$ . Then there exist $t_1,t_2$ such
that $t_1< t_2$ and $\phi_-(t-\sigma_0) < \phi_+(t), \ t \in
\mathcal{I}:= (-\infty, t_1] \cup [t_2, +\infty)$. Now, set $\sigma
= \sigma_0 + (t_2-t_1)$. Since both functions are increasing, we
have
$$
\phi_-(t- \sigma) \leq \phi_-(t- \sigma_0) < \phi_+(t),\quad t \in
\mathcal{I},
$$
$$
\phi_-(t- \sigma) < \phi_+(t-(t_2-t_1) ) \leq \phi_+(t), \quad t \in
[t_1,t_2].\qquad\square
$$
\end{pf} \vspace{-10mm}
\section{Proof of Theorem \ref{main}}\vspace{-3mm} \label{sm}
{\bf 6.1. Necessity.}  Let $u(t,x)= \zeta(ct+ \nu\cdot x)$ be a
positive bounded monotone solution of the delayed KPP-Fisher
equation. Then $\varphi(t)=\zeta(ct)$ satisfies \vspace{-3mm}
\begin{eqnarray} \label{twe2a}
&&  \epsilon \varphi''(t) - \varphi'(t) + \varphi(t)(1-
\varphi(t-h)) =0, \quad t \in \R,
\\ \nonumber
&& \epsilon \varphi'(t)= \epsilon \varphi'(0)-\varphi(0) +
\varphi(t) + \int_0^t\varphi(s)(1- \varphi(s-h))ds.
\end{eqnarray}
The latter relation implies that $\varphi(\pm\infty) \in \{0,1\}$
since otherwise $\varphi'(\pm\infty)=\infty$. Hence $\varphi:\R \to
(0,1)$. Let  $\phi \in C^2(\R,(0,1))$ be an arbitrary solution of
(\ref{twe2a}). Suppose for a moment that $\phi'(t_0)=0$. Then
necessarily $\phi''(t_0)<0$ so that $t_0$ is the unique critical
point (absolute maximum) of $\phi$. But then $\phi'(s) < 0$ for
 $s> t_0$, so that $\phi''(s) < 0, \ s \geq t_0$, which yields the
contradiction $\phi(+\infty)=-\infty$. In consequence, either
$\phi'(s)
>0$ or $\phi'(s) < 0$ for all $s \in \R$. But as we have seen, $\phi'(s) < 0$
implies $\phi(+\infty)=-\infty$, a contradiction. Hence, any
solution $\phi \in C^2(\R,(0,1))$ of (\ref{twe2a}) satisfies
$\phi'(t)>0,\ \phi(-\infty)=0, \ \phi(+\infty)=1$.
\begin{lem}\label{mon} If $\phi\in C^2(\R, (0,1))$ satisfies
(\ref{twe2a}), then $\epsilon \in (0,0.25]$.
\end{lem}\vspace{-5mm}\begin{pf} Suppose for a moment that $\epsilon > 0.25$. Then
the characteristic equation $\epsilon \lambda^2-\lambda+1=0$
associated with the trivial steady state of  (\ref{twe2a}) has two
simple complex conjugate roots $\omega_{\pm}= (2\epsilon)^{-1}(1\pm
i \sqrt{4\epsilon-1})$.

Since $\phi\in C^2(\R, (0,1))$ is a solution of (\ref{twe2a}), it
holds that $\phi'(t)>0, \ t \in \R, $ $\ \phi(-\infty)=0$. Set $z(t)
= (\phi(t),\phi'(t))^T$, it is easy to check that $z(t)$ satisfies
the following asymptotically autonomous linear differential equation
$$
z'(t) = (A+R(t))z(t), \ t \in \R, \quad A= \left(%
\begin{array}{cc}
  0 & 1 \\
  -1/\epsilon & 1/\epsilon \\
\end{array}%
\right), \ R(t) = \left(%
\begin{array}{cc}
0 & 0 \\
 \phi(t-h)/\epsilon & 0 \\
\end{array}%
\right).
$$
Since $R(-\infty)=0, \ \int_{-\infty}^0|R'(t)|dt = \phi(-h)$ and the
eigenvalues $\omega_{\pm}$ of $A$ are complex conjugate, we can
apply the Levinson theorem \cite[Theorem 1.8.3]{east} to obtain the
following asymptotic formulas at $t=-\infty$: \vspace{-5mm}
\begin{eqnarray*}
\phi(t) = (a+ o(1))e^{t/(2\epsilon)}
\cos(t\sqrt{4\epsilon-1}(1+o(1))+b +o(1)), \\
\phi'(t) = (c+ o(1))e^{t/(2\epsilon)}
\sin(t\sqrt{4\epsilon-1}(1+o(1))+d +o(1)),
\end{eqnarray*}
where $a^2+c^2 \not=0$. But this means that either $\phi(t)$ or
$\phi'(t)$ is oscillating around zero, a contradiction. \hfill
{$\square$}\end{pf}\vspace{-3mm}
\begin{lem}\label{c1b2} If $h > h_1$ or $h \in (1/e,h_1]$ and
$c> c^*(h)$ then Eq. (\ref{twe2a}) does not have any solution
$\phi\in C^2(\R, (0,1))$.
\end{lem}\vspace{-7mm}
\begin{pf} On the contrary, let us assume that Eq. (\ref{twe2a}) has a solution $\phi\in
C^2(\R, (0,1))$. Then Lemma \ref{mon} implies that $\epsilon \in
(0,0.25]$ and therefore the assumptions of this lemma imply that
$\psi(z,\epsilon)$ does not have negative zeros.  Following the
approach in \cite{TAT}, we will show that this will force $\phi(t)$
to oscillate about the positive equilibrium. For the convenience of
the reader, the proof is divided in several steps.

{\it Claim I:  $y(t):=1-\phi(t)  > 0$ has at least  exponential
decay as $t \to +\infty$.} First, observe that
\begin{equation}\label{ph} \epsilon y''(t) - y'(t)= \phi(t)y(t-h),
\quad t \in \R.
\end{equation}\vspace{-4mm}
Therefore, with $\gamma := \phi(t_0)$, which is close to $1$, and
$g(t):= \phi(t)y(t-h)- \phi(t_0)y(t)$, we obtain that
$$ \epsilon y''(t) - y'(t) - \gamma y(t) - g(t) =0, \
t \in \R.$$ Note  that  $g(t) > 0$ for all sufficiently large $t$.
Since $y(t), g(t)$ are bounded on $\R$, it holds that
$$
y(t)= -\frac{1}{\epsilon(m-l)}\left(\
\int_{-\infty}^te^{l(t-s)}g(s)ds +
\int_t^{+\infty}e^{m(t-s)}g(s)ds\right),
$$
where $l<0$ and $0<m$ are roots of  $\epsilon z^2 - z - \gamma =0$.
The latter representation of $y(t)$ implies that there exists $T_0$
such that\vspace{0mm}
\begin{equation}\label{yyp}
y'(t)- ly(t) = - \frac{1}{\epsilon} \int_t^{+\infty}e^{m(t-s)}g(s)ds
< 0, \ t \geq T_0.
\end{equation}\vspace{-3mm}
Hence, $(y(t)\exp(-lt))' < 0,\  t \geq T_0,$ and therefore
\begin{equation}\label{yg}
y(t) \leq y(s)e^{l(t-s)}, \quad t\geq s \geq T_0, \quad g(t)=
O(e^{lt}), \ t \to +\infty.
\end{equation}
\vspace{-3mm} It is easy to see that these estimates are valid for
every negative $l> (2\epsilon)^{-1}(1-\sqrt{1+4\epsilon})$. Finally,
(\ref{yyp}), (\ref{yg}) imply that $y'(t)= O(e^{lt}), \ t \to
+\infty$.

{\it Claim II: $y(t):=1-\phi(t)  > 0$ is not
superexponentially small as $t \to +\infty$.}\\
We already have proved that $y(t)$ is strictly decreasing and
positive on  $\R$. Since the right hand side of Eq. (\ref{ph}) is
positive and integrable on $\R_+$, and since $y(t)$ is a bounded
solution of (\ref{ph}) satisfying $y(+\infty)=0$, we find that
\begin{equation}\label{2r}
y(t) = \int_t^{+\infty}(1-e^{(t-s)/\epsilon})\phi(s)y(s-h)ds.
\end{equation}\vspace{-5mm}
As a consequence, there exists $T_1$ such that
$$
y(t) \geq 0.5(1-e^{-0.5h/\epsilon}) \int_{t-0.5h}^{t}y(s)ds:= \xi
\int_{t-0.5h}^{t}y(s)ds, \quad t \geq T_1-h.
$$
Now, since $y(t) >0$ for all $t$, we can find positive $C, \rho$
such that $y(s) > C e^{-\rho s}$ for all $s \in [T_1-h, T_1]$. We
can assume that $\rho$ is large enough to satisfy the inequality $
\xi(e^{0.5\rho h}-1)> {\rho}$. Then we claim that $y(s) > C e^{-\rho
s}$ for all $s \geq T_1-h$. Conversely, suppose that $t'
> T_1$ is the leftmost point where $y(t') = C e^{-\rho t'}$. Then we get a contradiction:
$$
y(t') \geq \xi \int_{t'-0.5h}^{t'}y(s)ds > C\xi
\int_{t'-0.5h}^{t'}e^{-\rho s}ds = Ce^{-\rho t'}\xi\frac{e^{0.5\rho
h}-1}{\rho} > Ce^{-\rho t'}.$$ {\it Claim III: $y(t)> 0$ can not
hold when $\psi(z,\epsilon)$ does not have any zero in
$(-\infty,0)$}. Observe that $y(t)=1- \phi(t)$ satisfies
$$\epsilon y''(t) - y'(t)- (1-y(t))y(t-h)=0,\ t \in \R,\
$$
where in virtue of Claim I, it holds that $(y(t),y'(t))= O (lt)$ at
$t = +\infty$.  Then \cite[Proposition 7.2]{FA}  implies that there
exists $\gamma <l$ such that $y(t) = v(t) + O(\exp(\gamma t)),$ $ t
\to +\infty,$ where $v$ is a {\it non empty} (due to Claim II)
finite sum of eigensolutions of the limiting equation
$$
\epsilon y''(t) - y'(t)- y(t-h)=0,\ t \in \R,\
$$
associated to the eigenvalues $\lambda_j \in F= \{\gamma <
\Re\lambda_j \leq l\}$. Now, since the set $F$ does not contain any
real eigenvalue by our assumption, we conclude that $y(t)$ should be
oscillating on $\mathbb{R}_+$,  a contradiction. \hfill
{$\square$}\end{pf} \vspace{-5mm}

{\bf 6.2. Sufficiency.} Suppose that $\epsilon \in (0, 0.25]$ and
let $0 < \lambda \leq  \mu$ be the roots of the equation $\epsilon
z^2 -z + 1 =0$. In Lemmas \ref{c1b}-\ref{c1d} below, $\mathcal{K}$
stands either for $\mathcal{A}$ or $\mathcal{B}$ (defined by
(\ref{iie}), (\ref{iie2})).
\begin{lem}\label{c1b} If $\phi, \psi \in C(\R, (0,1))$ and
$\phi(t) \leq  \psi(t)$ for all $t \in \R$, then $\mathcal{K}\phi,
\mathcal{K}\psi \in C(\R, (0,1))$ and $(\mathcal{K}\phi)(t) \leq
(\mathcal{K}\psi)(t), \ t \in \R$. Moreover, if $\phi$ is increasing
then $\mathcal{K}\phi$ is also increasing.
\end{lem}
\vspace{-7mm}
\begin{pf} The proof is straightforward. \hfill {$\square$}\end{pf}
\begin{lem}\label{c1c} Let $\epsilon \in (0,0.25]$. If $\phi_+ \in C^1(\R, (0,1))$ satisfies
the inequality
$$
\epsilon \phi''(t) - \phi'(t) + \phi(t)(1- \phi(t-h)) \geq 0
$$
for all $t \in \R':=\R\setminus\{T_1,\dots,T_m\}$ and $\phi_+''(t),
\phi_+'(t)$ are bounded on $\R'$,  then $(\mathcal{K}\phi_+)(t) \leq
\phi_+(t)$ for all $t \in \R$.
\end{lem}\vspace{-5mm}\begin{pf} If $\omega(T_i):=0$ and
$$
\omega(t): = \epsilon \phi_+''(t) - \phi_+'(t) + \phi_+(t)(1-
\phi_+(t-h)), \quad t \in \R'=\R\setminus\{T_1,\dots,T_m\}
$$
then $\omega(t)\geq 0$ for all $t \in \R',$ $\omega(t)$ is bounded
on $\R'$ and
$$
\epsilon \phi_+''(t) - \phi_+'(t) + \phi_+(t) = \omega_1(t), \quad t
\in \R',
$$
where $\omega_1(t): =\omega(t)+ \phi_+(t)\phi_+(t-h)$ is bounded on
$\R'$. Let now $\epsilon \in (0,0.25)$. By Lemma \ref{imp}, we
obtain that \vspace{-5mm}
\begin{eqnarray*}
  \phi_+(t) &=& \frac{1}{\epsilon(\mu - \lambda)}
\int_t^{+\infty}\left(e^{\lambda (t-s)}- e^{\mu
(t-s)}\right)\omega_1(s)ds =  \\
   & & (\mathcal{A}\phi_+)(t)+
\frac{1}{\epsilon(\mu - \lambda)} \int_t^{+\infty}\left(e^{\lambda
(t-s)}- e^{\mu (t-s)}\right)\omega(s)ds \geq (\mathcal{A}\phi_+)(t).
\end{eqnarray*}
The case $\epsilon=0.25$ (which corresponds to
$\mathcal{K}=\mathcal{B}$) is completely analogous to the previous
one. \hfill {$\square$}\end{pf}\vspace{-5mm} The proof of the next
lemma is similar to that of Lemma \ref{c1c}:
\begin{lem}\label{c1f} Let $\epsilon \in (0,0.25]$. If $\phi_- \in C^1(\R, (0,1))$ satisfies
the inequality
$$
\epsilon \phi''(t) - \phi'(t) + \phi(t)(1- \phi(t-h)) \leq 0
$$
for all $t \in \R\setminus\{T_1,\dots,T_m\}$ and $\phi_-''(t),
\phi_-'(t)$ are bounded on $\R\setminus\{T_1,\dots,T_m\}$,  then
$(\mathcal{K}\phi_-)(t) \geq \phi_-(t)$ for all $t \in \R$.
\end{lem}
Set $\phi_{j+1}^{\pm}:= (\mathcal{K}\phi_{j}^{\pm}), \ j \geq 0, \
\phi_{0}^{\pm}: = \phi_{\pm},$ and let the increasing functions
$\phi_{-}\leq \phi_{+}$ be as in Lemmas \ref{c1c}, \ref{c1f}. Then
$$ \phi_{-} \leq \phi_{1}^{-}\leq \dots \leq \Phi_- \leq
\Phi_+\leq\dots \phi_{j}^{-} \dots \leq
  \dots \phi_{1}^{+} \leq \phi_{+},
$$
where $\Phi_{\pm}(t) = \lim_{j\to \infty} \phi_{j}^{\pm}(t)$
pointwise and $\phi_{j}^{\pm}$ are increasing (by Lemma \ref{c1b}).
\begin{lem}\label{c1d} $\Phi_\pm$ are wavefronts and $\Phi_{\pm}(t) = \lim_{j\to \infty} \phi_{j}^{\pm}(t)$
uniformly on $\R$.
\end{lem}\vspace{-4mm}
\begin{pf}
Applying the Lebesgue's dominated convergence theorem to
$\phi_{j+1}^-:= \mathcal{K}\phi_{j}^-$, we obtain that $\Phi_{-}(t)=
(\mathcal{K}\Phi_{-})(t)$. Differentiating this equation twice with
respect to $t$ , we deduce that $\Phi_{-}:\R\to (0,1)$ is a
$C^2$-solution of (\ref{twe2a}) (and thus $\Phi_{-}'(t)
>0$). As a consequence of the Dini's theorem, we have that
$\Phi_-(t) = \lim_{j\to \infty} \phi_{j}^-(t)$ uniformly on compact
sets. Since $\Phi_-, \phi_{j}^-$ are asymptotically constant and
increasing, this convergence is uniform on $\R$. The proof for
$\Phi_+$ is similar. \hfill {$\square$}\end{pf} \vspace{-7mm}
\begin{cor} \label{mcor} Eq. (\ref{17}) has a  monotone  wavefront
$u(x,t) = \zeta(x\cdot\nu + ct),$ \ $ |\nu|=1,$ connecting $0$ with
$1$ if one of the following conditions holds
\begin{enumerate}
  \item $0\leq h \leq 1/e$ and $2 \leq c$;
  \item $1/e < h <  \ h_1$ and $2 \leq c < c^*(h)$.
\end{enumerate}
\end{cor}\vspace{-5mm}
\begin{pf}  It is an immediate consequence of Lemmas \ref{low},
\ref{up}, \ref{as}, \ref{c1c}-\ref{c1d}. \hfill {$\square$}\end{pf}
\vspace{-5mm} If $c=c^*(h)$,  the reasoning of the last proof does
not apply because of the lack of explicit upper solutions. Below, we
follow an idea from \cite[Section 6]{TAT}:
\begin{lem}\label{mlem} Eq. (\ref{17}) has a positive monotone  wavefront
$u(x,t) = \zeta(x\cdot\nu + ct),$ $|\nu|=1,$ connecting $0$ with $1$
if  $1/e < h \leq   \ h_1$ and $c = c^*(h)$.
\end{lem}
\vspace{-5mm}
\begin{pf}
\underline{Case I.} Fix some $h \in (1/e,h_1)$ and $\epsilon=
\epsilon^*(h)$. Then there exists a decreasing sequence $\epsilon_j
\downarrow \epsilon^*(h)$  such that Eq. (\ref{twe2a}) has at least
one monotone positive heteroclinic solution $\phi_j(t)$ normalized
by $\phi_j(0)=0.5$. It is clear that
$\phi_j(t)=(\mathcal{A}\phi_j)(t)$. Moreover, each
$y_j(t):=1-\phi_j(t)  > 0$ solves (\ref{2r}) so that
$$
|\phi_j'(t)| =
|\frac{1}{\epsilon}\int_t^{+\infty}e^{(t-s)/\epsilon}\phi(s)(1-\phi_j(s-h))ds|\leq
1, \ t \in \R. $$ Thus, by the Ascoli-Arzel${\rm \grave{a}}$ theorem
combined with the diagonal method, $\{\phi_j\}$ has a subsequence
$\{\phi_{j_k}\}$ converging (uniformly on compact subsets of $\R$)
to some continuous non-decreasing non-negative function $\phi_*, \
\phi_*(0) = 0.5$. Applying the Lebesgue's dominated convergence
theorem to $\phi_{j_k}(t)=(\mathcal{A}\phi_{j_k})(t)$, we find that
$\phi_*$ is also a fixed point of $\mathcal{A}$. Hence, $\phi_*: \R
\to [0,1]$ is a monotone solution of Eq. (\ref{twe2a}) considered
with $\epsilon= \epsilon^*(h)$. Since $\phi_*(0) = 0.5$, $\phi_*: \R
\to (0,1)$ is actually a monotone wavefront.

\underline{Case II.} Finally, let $\epsilon= 0.25$ and $h=h_1$. This
case can be handled exactly in the same way as Case I if we keep
$\epsilon= 0.25$ fixed, replace $\mathcal{A}$ with $\mathcal{B}$,
and take some increasing sequence $h_j \uparrow h_1$ instead of
$\epsilon_j \downarrow \epsilon^*(h)$. \hfill {$\square$}\end{pf}
\begin{cor} \label{mcoro} Assume that $c = c^*(h), \ 1/e < h \leq  h_1,$
and let $\phi_{-}$ be as in (\ref{sokr}). If $A$ is sufficiently
large, then
$$
\phi_{-} \leq \phi_{1}^{-}\leq \dots \leq \phi_{j}^{-} \dots \leq
\Phi = \mathcal{K}\Phi,
$$
where $\Phi$ is a wavefront and $\Phi(t) = \lim_{j\to \infty}
\phi_{j}^-(t)$ uniformly on $\R$.
\end{cor}
\begin{pf}
If $c = c^*(h)$ we will take the heteroclinic solution $\Phi_-$
whose existence was established in Lemma \ref{mlem} as an upper
solution. Due to (\ref{3up}), we can assume  that
$$\beta_{+,1}:= \lim_{t \to +\infty}(1-\Phi_{-}(t))t^{-1}e^{-\lambda_2
t} >0.$$ Next, let $\phi_-$ be defined by (\ref{sokr}). Since
$\alpha_- := \lim_{t \to -\infty}\phi_{-}(t)e^{-\lambda t} = 0$ and
$$\beta_{-,1}: = \lim_{t \to
+\infty}(1-\phi_{-}(t-\frac{1}{A}))t^{-1}e^{-\lambda_2 t}
=Ae^{-\lambda_2/A} > \beta_{+,1},$$ for sufficiently large $A$,
Lemma \ref{as} implies that $\phi_-(t) < \Phi_-(t+\sigma), \ t \in
\R,$ for some $\sigma$. Finally, it suffices to take $\phi_+(t):=
\Phi_-(t+\sigma)$ and repeat the proof of Lemma \ref{c1d}. \hfill
{$\square$}\end{pf} \vspace{-3mm} {\bf 6.3. Uniqueness.} Our method
of proof follows a nice idea due to Diekmann and Kaper, see
\cite[Theorem 6.4]{DK}.  Suppose that $c\not= c^*(h)$ and let
$\phi_1,\phi_2$ be two different (modulo translation) profiles of
wavefronts propagating at the same speed $c$. Due to Theorem
\ref{mainas}, we may assume that $\phi_1, \phi_2$ have  the same
asymptotic representation $\phi_j(t)= 1-e^{\lambda_2t}(1+o(1))$ at
$+\infty$. Moreover, $\phi_j= \mathcal{K}\phi_j$, where
$\mathcal{K}= \mathcal{A}$ if $c
>2$ and $\mathcal{K}= \mathcal{B}$ if $c =2$. Set $\omega(t):=
|\phi_2(t)-\phi_1(t)|e^{-\lambda_2t}$. Then $\omega(\pm\infty)=0,\
\omega(t)\geq 0, \ t \in \R,$ and $\omega(\tau)= \max_{s \in \R}
\omega(s): = |\omega|_0 >0$  for some $\tau$. From the identity
$\phi_2-\phi_1= \mathcal{K}\phi_2-\mathcal{K}\phi_1$, we deduce that
$$
  \omega(\tau) < \frac{e^{-\lambda_2\tau}}{\epsilon(\mu-\lambda)}
\int_\tau^{+\infty}(e^{\lambda (\tau-s)}- e^{\mu
(\tau-s)})(\omega(s)e^{\lambda_2s}+\omega(s-h)e^{\lambda_2(s-h)})ds
<$$
$$
 \frac{|\omega|_0e^{-\lambda_2\tau}}{\epsilon(\mu-\lambda)}
\int_\tau^{+\infty}(e^{\lambda (\tau-s)}- e^{\mu
(\tau-s)})(e^{\lambda_2s}+e^{\lambda_2(s-h)})ds = |\omega|_0=
\omega(\tau), \ {\rm if }\ c
>2; $$
$$
  \omega(\tau) < 4e^{-\lambda_2\tau}
\int_\tau^{+\infty}(s-\tau)e^{\lambda
(\tau-s)}(\omega(s)e^{\lambda_2s}+\omega(s-h)e^{\lambda_2(s-h)})ds
<$$
$$
 4|\omega|_0e^{-\lambda_2\tau}
\int_\tau^{+\infty}(s-\tau)e^{\lambda
(\tau-s)}(e^{\lambda_2s}+e^{\lambda_2(s-h)})ds = |\omega|_0=
\omega(\tau), \ \ {\rm if }\ c =2, $$ which is impossible. Hence,
$|\omega|_0=0$ and the proof is complete. \hfill {$\square$}
\section{Proof of Theorem \ref{mainas}}\label{rrra}
First, using the bilateral Laplace transform
$(\mathcal{L}y)(z):=\int_{\R}e^{-sz}y(s)ds$ (see e.g.
\cite{widder}), we extend \cite[Proposition 7.1]{FA} (see also
\cite[Lemma 4.1]{AVT} and \cite[Lemma 22]{TAT}) for the case $J =
\R$.
\begin{lem}\label{41lin}  Set $\chi(z):=
z^2+\alpha z+\beta +pe^{-zh}$ and let $y\in C^2(\R,\R)$ satisfy
\begin{equation}\label{twel}
y''(t) +\alpha y'(t)+\beta y(t)+ py(t-h) = f(t), \ t \in \R,
\end{equation}
where  $\alpha,\beta,p,h \in \R$ and
\begin{equation}\label{2u}
y(t)=\left\{%
\begin{array}{ll}
     O(e^{-Bt}),
    & \hbox{as}\ t \to +\infty,   \\
     O(e^{bt}), & \hbox{as}\ t \to -\infty;
\end{array}%
\right.
f(t)=\left\{%
\begin{array}{ll}
     O(e^{-Ct}),
    & \hbox{as}\ t \to +\infty,   \\
     O(e^{ct}), & \hbox{as}\ t \to -\infty,
\end{array}%
\right.
\end{equation}
for some non-negative $b < c, B < C, b+B >0$.   Then, for each
sufficiently small $\sigma >0$, it holds that
$$
y(t)=\left\{%
\begin{array}{ll}
     w_+(t) +e^{-(C-\sigma)t}o(1),
    & \hbox{as}\ t \to +\infty,    \\
     w_-(t) + e^{(c-\sigma)t}o(1), & \hbox{as}\ t \to -\infty,
\end{array}%
\right.
$$
where
$$w_{\pm}(t)=\pm \sum_{ \lambda_j \in F_\pm} {\rm Res}_{z=\lambda_j}
\left[\frac{e^{zt}}{\chi(z)}\int_{\R}e^{-zs}f(s)ds\right]
$$ is a finite sum of  eigensolutions of equation (\ref{twel}) associated to
the eigenvalues $\lambda_j \in F_+ = \{-C+\sigma < \Re\lambda_i \leq
-B\}$ and $\lambda_j \in F_- = \{b \leq \Re\lambda_i <c-\sigma\}$.
\end{lem}\begin{pf} We will divide our proof into several parts.

\underline{Step I.} We claim that there exist non-negative $B',b'$
such that $B' \leq B, b'\leq b,$ $B'+b' >0$ and
\begin{equation}\label{2uu}
y'(t), y''(t) =\left\{%
\begin{array}{ll}
     O(te^{-B't}),
    & \hbox{as}\ t \to +\infty,   \\
     O(te^{b't}), & \hbox{as}\ t \to -\infty.
\end{array}%
\right.
\end{equation}
We will distinguish two cases:

\underline{Case A}. Suppose that $\alpha=0$. Then clearly
$y''(t)=O(e^{-Bt})$ at $t= -\infty$, is bounded on $\R$ and
therefore $y'(t)$ is uniformly continuous on $\R$. Since $B+b>0$
then either $y(+\infty)=0,\ \limsup_{s \to -\infty}|y(s)| < \infty$
or $y(-\infty)=0,\ \limsup_{s \to +\infty}|y(s)| < \infty$. Suppose,
for example that $B>0$ (hence $y(+\infty)=0$), the other case being
similar. Then, applying the Barbalat lemma, see  e.g. \cite{wz}, we
find that $y'(+\infty)=0$. This implies that $y'(t)= -
\int_t^{+\infty}y''(s)ds= O(e^{-Bt})$ at $t= +\infty$. Thus we may
set $B'=B$. Now, $y'(t)= y'(0)+ \int_0^{t}y''(s)ds= O(t)$ at $t=
-\infty$ so that we can choose $b' =0$.

\underline{Case B}. Let now $\alpha\not=0$. For example, suppose
that $\alpha >0$ (the case $\alpha <0$ is similar).  Then, for some
$\xi$,
$$
y'(t)= \xi e^{-\alpha t} + \int_{-\infty}^te^{-\alpha (t-s)}
\{f(s)-\beta y(s)- py(s-h)\}ds.
$$
In fact, since the second term of the above formula is bounded on
$\R$ and we can not have $y'(-\infty)= \pm \infty$ (due to the
boundedness of $y(t)$), we obtain that $\xi=0$. But then $y'(t)=
O(e^{bt}),\ t \to -\infty$ and $y'(t)= O(te^{-\min\{\alpha, B\}t}),\
t \to +\infty$. Note that $b'+B' = \min\{\alpha+b, B+b\}
>0$. Finally, (\ref{twel}) assures that (\ref{2uu}) is also valid
for $y''(t)$.

\underline{Step II.} Applying the bilateral Laplace transform
$\mathcal{L} $ to (\ref{twel}), we obtain that $
\chi(z)\tilde{y}(z)= \tilde{f}(z), $ where $\tilde{y} =
\mathcal{L}y, \ \tilde{f} = \mathcal{L}f$ and $-B' < \Re z  < b'$.
Moreover, from the growth restrictions (\ref{2u}), we conclude that
$\tilde{y}$ is analytic in $-B<\Re z < b$ while $\tilde{f}$ is
analytic in $-C <\Re z <c$. As a consequence, $H(z)=\tilde{f}(z)
/\chi(z)$ is analytic in $-B<\Re z < b$ and meromorphic in $-C <\Re
z <c$. Observe that $H(z) = O(z^{-2}),\ z \to \infty,$ for each
fixed strip $\Pi(s_1,s_2)=\{s_1 \leq \Re z \leq s_2\},$ $ -C <
s_1<s_2<c$. Now, let $ \sigma > 0$ be such that the vertical strips
$c - 2\sigma < \Re z < c$ and $-C < \Re z < -C+ 2\sigma$ do not
contain any zero of $\chi(z)$. By the inversion formula
\cite[Theorem 5a]{widder}, for each $\delta \in (-B,b)$, we obtain
that
$$
y(t) = \frac{1}{2\pi i}\int_{\delta -i\infty}^{\delta
+i\infty}e^{zt}\tilde{y}(z)dz = \frac{1}{2\pi i}\int_{\delta
-i\infty}^{\delta +i\infty}e^{zt}H(z)dz =w_{\pm}(t) + u_\pm(t),  \ t
\in \R,
$$
$$
\hspace{-3mm}{\rm where} \ w_\pm(t) = \pm \sum_{ \lambda_j \in
F_\pm} {\rm Res}_{z=\lambda_j} \frac{e^{zt}\tilde{f}(z) }{\chi(z)},
\ u_\pm(t) = \frac{1}{2\pi i}\int_{\mp(c -\sigma) -i\infty}^{\mp(c
-\sigma) +i\infty}e^{zt}H(z)dz.
$$
The above sum is finite, since $\chi(z)$ has a finite set of the
zeros in $F_{\pm}$. Now, for $a(s) = H(\mp(c-\sigma) + is)$, we
obtain that
$$
u_{\pm}(t) =\frac{e^{\mp(c - \sigma)t}}{2\pi}\left\{
\int_{\R}e^{ist}a(s)ds \right\}, \ t \in \R.$$ Next, since $a \in
L_1(\R)$, we have, by the Riemann-Lebesgue lemma, that
$$
\lim_{t \to \infty}\int_{\R}e^{ist}a_1(s)ds =0.
$$
Thus we get $u_{\pm}(t) = e^{\mp(c - \sigma)t} o(1)$ at $t =
\infty$, and the proof is completed. \hfill {$\square$}\end{pf}

Now we can prove Theorem \ref{mainas}:
\begin{pf}[Theorem \ref{mainas}] \underline{Case I: asymptotics at $t=
+\infty$}. It follows from (\ref{yg}) that $y(t)=1- \phi(t)$
satisfies $y(t) = O(e^{lt}), \ t \to +\infty,$ for every negative $l
> (2\epsilon)^{-1}(1-\sqrt{1+4\epsilon})$. Moreover, $f(t):=
-y(t)y(t-h) = O(e^{2lt}),  t \to +\infty, $ $\ y(t) = O(1), \ t \to
-\infty$ and
$$\epsilon y''(t) - y'(t)- y(t-h)=-y(t)y(t-h),\ t \in \R.\
$$
Therefore Lemma \ref{41lin} implies that, for every small
$\sigma>0$,
$$
y(t)=\sum_{ 2l+\sigma < \Re\lambda_j \leq l} {\rm Res}_{z=\lambda_j}
\frac{e^{zt}\tilde{f}(z) }{\chi(z)}+ e^{(2l+\sigma)t}o(1), \ t \to
+\infty.
$$
Now, observe that $(2\epsilon)^{-1}(1-\sqrt{1+4\epsilon}) >
\lambda_2$ so that either $\lambda_2 \in (2l+\sigma,l)$ or
$\lambda_2 \leq 2l$. In the latter case, we obtain
$y(t)=e^{(2l+\sigma)t}o(1), \ t \to +\infty,$ which allows to repeat
the above procedure till the inclusion $\lambda_2 \in
(2^jl+\sigma,2^{j-1}l)$ is reached for some integer $j$. In this
way, assuming  that $\lambda_1 < \lambda_2$, for each small $\sigma
>0$, we find that
\begin{equation}\label{for1}\hspace{-7mm}
y(t)= \eta e^{\lambda_2t} +  O(e^{(\lambda_2-\sigma)t}),  \ {\rm
where \ }\eta:=
\frac{\int_{\R}e^{-\lambda_2s}y(s)y(s-h)ds}{-\chi'(\lambda_2)}
>0.
\end{equation}
Now, if $c = c^*(h)$ (i.e. $\lambda_1 = \lambda_2$), we obtain
analogously that
$$
y(t+t_0)= \xi te^{\lambda_2t} +  O(e^{(\lambda_2-\sigma)t}), \ t \to
+\infty,
$$
for some appropriate $t_0$ and $\xi >0$.

Suppose now that $h  \in (0,h_0], \ c\leq c^\#(h)$. Then Lemmas
\ref{c1dd}, \ref{c1de} imply that $\Re \lambda_j < \lambda_1 \leq
2\lambda_2$. This means that formula (\ref{for1}) can be improved as
follows:
$$
y(t)= \eta e^{\lambda_2t} + O(e^{(2\lambda_2+\sigma)t}), \ t \to
+\infty.
$$
Finally, if  $h \in (0.5\ln2, h_0]$ and $c \in (c^\#(h), c^*(h))$,
it holds that $2\lambda_2 < \lambda_1 < \lambda_2$. Then
\begin{eqnarray*}
  y(t) &=& \sum_{ 2\lambda_2+\sigma < \Re\lambda_j \leq \lambda_2} {\rm
Res}_{z=\lambda_j} \frac{e^{zt}\tilde{f}(z) }{\chi(z)}+
e^{(2\lambda_2+\sigma)t}o(1) =   \\
  & & \hspace{-17mm} \eta e^{\lambda_2t} + \theta e^{\lambda_1t}+ e^{(\lambda_1-\sigma)t}o(1), \ {\rm where \ }\theta:=
\frac{\int_{\R}e^{-\lambda_1s}y(s)y(s-h)ds}{-\chi'(\lambda_1)} <0.
\end{eqnarray*}
\underline{Case II: asymptotics at $t= -\infty$}. This case is much
easier to analyze since  the characteristic polynomial $\epsilon z^2
- z+1$ of the variational equation
\begin{equation}\label{veq}
\epsilon y''(t) - y'(t) + y(t)=0, \ \epsilon \in (0,0.25],
\end{equation}
along the trivial equilibrium of (\ref{twe2a}) has only two real
zeros $0 < \lambda \leq \mu$. It is easy to check that $2\lambda
\leq \mu$ if and only if $c\geq 1.5\sqrt{2}=2.121\dots$.

Since $\phi(-\infty)=0$ and equation (\ref{veq}) is exponentially
unstable on $\R_-$, we conclude that the perturbed equation
$$
\epsilon y''(t) - y'(t) + y(t)(1-\phi(t-h)) =0
$$
is also exponentially unstable on $\R_-$ (e.g. see \cite{dal}). As a
consequence, $\phi(t) = O(e^{mt}),\ t \to -\infty,$ for some $m>0$.
Now we can proceed as in Case I, since
$$
\epsilon \phi''(t) - \phi'(t) + \phi(t) = f_1(t),
$$
with $f_1(t):=\phi(t)\phi(t-h) = O(e^{2mt})$.  The details are left
to the reader. \hfill {$\square$}\end{pf}
\section*{Acknowledgments}
The authors thank Teresa Faria,  Anatoli Ivanov and Eduardo Liz for
useful discussions.  Sergei Trofimchuk was partially supported by
CONICYT (Chile) through PBCT program ACT-05 and by the University of
Talca, through program ``Reticulados y Ecuaciones". Research was
supported in part by FONDECYT (Chile), project 1071053.

\end{document}